\documentclass[11pt]{amsart}

\usepackage{amsmath}
\usepackage{amssymb}

\usepackage{amsfonts}
\usepackage{amsthm}
\usepackage{enumerate}
\usepackage[mathscr]{eucal}

\usepackage{tikz}

\usepackage{bbm}
\usepackage[caption=false,font=footnotesize]{subfig}

\usepackage{graphicx}
\usepackage{verbatim}

\newcommand{\Be}{\begin{equation}}
\newcommand{\Ee}{\end{equation}}

\newcommand{\Bea}{\begin{eqnarray}}
\newcommand{\Eea}{\end{eqnarray}}
\newcommand{\Beas}{\begin{eqnarray*}}
\newcommand{\Eeas}{\end{eqnarray*}}
\newcommand{\Benu}{\begin{enumerate}}
\newcommand{\Eenu}{\end{enumerate}}
\newcommand{\Bi}{\begin{itemize}}
\newcommand{\Ei}{\end{itemize}}

\def\intslash{\rlap{\kern  .32em $\mspace {.5mu}\backslash$ }\int}
\def\qsl{{\rlap{\kern  .32em $\mspace {.5mu}\backslash$ }\int_{Q_x}}}

\def\rr{\mathbb R}

\def\emph#1{{\it #1 }}

\def\supp{{\text{\rm supp}}}

\def\inn#1#2{\langle#1,#2\rangle}

\def\lc{\lesssim}

\def\bbD{{\mathbb {D}}}
\def\bbE{{\mathbb {E}}}

\def\bbN{{\mathbb {N}}}

\def\bbR{{\mathbb {R}}}

\def\bbZ{{\mathbb {Z}}}

\def\f{{\overline{f}}}
 at 10 true pt

\def\be#1{\begin{equation}\label{ #1}}
\def\endeq{\end{equation}}
\def\endal{\end{align}}
\def\bas{\begin{align*}}
\def\eas{\end{align*}}
\def\bi{\begin{itemize}}
\def\ei{\end{itemize}}

\def\emph#1{{\it #1}}
\def\textbf#1{{\bf #1}}
\usepackage{bbm}
\def\bbone{{\mathbbm 1}}

\theoremstyle{plain}
  \newtheorem{theorem}{Theorem}[section]

   \newtheorem{proposition}[theorem]{Proposition}

\theoremstyle{remark}

\theoremstyle{definition}

\begin{document}

\title
[Dyadic averaging operators ]
{On uniform  boundedness of dyadic averaging operators in spaces of Hardy-Sobolev type}

\author[G. Garrig\'os \ \ \ A. Seeger \ \ \ T. Ullrich] {Gustavo Garrig\'os   \ \ \ \   Andreas Seeger \ \ \ \ 
Tino Ullrich}

\address{Gustavo Garrig\'os\\ 
Department of Mathematics\\
University of Murcia\\
30100 Espinardo\\Murcia, 
Spain
} 
\email{gustavo.garrigos@um.es}

\address{Andreas Seeger \\ Department of Mathematics \\ University of Wisconsin \\480 Lincoln Drive\\ Madison, WI,
53706, USA} \email{seeger@math.wisc.edu}

\address{Tino Ullrich\\
Hausdorff Center for Mathematics\\ Endenicher Allee 62\\
53115 Bonn, Germany} \email{tino.ullrich@hcm.uni-bonn.de}
\begin{abstract} We give an alternative proof of recent results by the authors on
uniform  boundedness of dyadic averaging operators 
in (quasi-)Banach spaces of Hardy-Sobolev and Triebel-Lizorkin type. 
{This result served as the main tool to establish Schauder basis properties of suitable enumerations of the univariate Haar system in the mentioned spaces.} 
%This result has been already proved in a previous work by the authors entitled
%``The Haar system as a Schauder basis in spaces of Hardy-Sobolev type''. It implies that a suitable enumeration of 
%the Haar system is a Schauder basis in these function spaces. 
%We use 
The rather elementary proof here is based on 
characterizations of the respective spaces in terms of
orthogonal compactly supported Daubechies wavelets.
% as the main tool. 
\end{abstract}

%\date{October 1, 2016}

%\date\today

\thanks{G.G. was supported in part   by grants  MTM2013-40945-P and MTM2014-57838-C2-1-P
from MINECO (Spain), and grant 19368/PI/14  from Fundaci\'on S\'eneca (Regi\'on
de Murcia, Spain). A.S. was supported in part by NSF grant 
DMS 1500162. T.U. was supported 
the DFG  Emmy-Noether program UL403/1-1}
\maketitle

%\address{}

%\begin{thanks}\end{thanks}
%\begin{abstract}
%\end{abstract}
%\centerline{\bf Preliminary version, March 4, 2015}

\section{Introduction} 
%In this paper we will give an alternative proof of a boundedness result \cite[Cor.\ 1.3]{gsu-basic} for 
Consider the dyadic averaging operators $\bbE_N$ on the real line given by  
\Be
\label{condexp}
\bbE_N f(x)= \sum_{\mu\in \bbZ} \bbone_{I_{N,\mu}} (x) \, 2^N \int_{I_{N,\mu}} f(t) dt\, 
\Ee
with $I_{N,\mu}= [2^{-N}\mu,2^{-N}(\mu+1))$. 
$\bbE_N f$  is  the conditional expectation of $f$ 
%operators $\bbE_N$
with respect to the $\sigma$-algebra generated by the dyadic intervals of length $2^{-N}$. 
The following theorem on uniform boundedness in 
Triebel-Lizorkin spaces $F^s_{p,q}$ 
was  proved by the authors in \cite{gsu-basic} {and serves as the main tool to establish that 
suitably regular enumerations of the Haar system form a Schauder basis for the spaces $F^s_{p,q}$ in the parameter ranges of the theorem. } 
Since the {uniform} boundedness result is interesting on its own we give an alternative proof based on wavelet theory {to make it accessible for 
a broader readership}.
%in the largest possible range of parameters $1/2<p<\infty$, 
%$0<q\leq \infty$, $1/p-1<s<\min\{1/p,1\}$.

\begin{theorem}\label{expthm} \cite{gsu-basic} Let $1/2<p<\infty$, $0<q\leq\infty$, and $1/p-1<s<\min\{1/p,1\}$. 
Then there is a constant $C:=C(p,q,s)>0$ such that for all $f\in F^s_{p,q}$
\begin{equation}\label{E_N}
    \sup\limits_{N\in \bbN}\|\bbE_Nf\|_{F^s_{p,q}} \leq C\|f\|_{F^s_{p,q}}.
\end{equation}
\end{theorem}

In \cite{gsu-basic}, this result served as the main  tool to establish that suitably regular enumerations of the Haar system form a Schauder basis for the spaces $F^s_{p,q}$ in the parameter ranges of the theorem, see
 \S\ref{Schauder}. The connection with  the Haar system is given via 
 %connection to that problem is given by 
 the martingale difference operators 
 \[\bbD_N=\bbE_{N+1}-\bbE_N\]
 which are the orthogonal projections to the spaces generated by Haar functions with fixed Haar frequency $2^N$. 
 
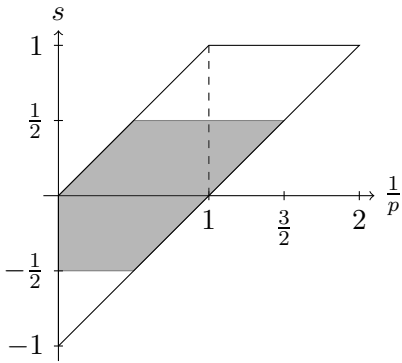
\begin{figure}[h]
 \begin{center}
\begin{tikzpicture}[scale=2]
% Axes
\draw[->] (-0.1,0.0) -- (2.1,0.0) node[right] {$\frac{1}{p}$};
\draw[->] (0.0,-1.1) -- (0.0,1.1) node[above] {$s$};

% Ticks
\draw (1.0,0.03) -- (1.0,-0.03) node [below] {$1$};
\draw (2.0,0.03) -- (2.0,-0.03) node [below] {$2$};
\draw (1.5,0.03) -- (1.5,-0.03) node [below] {$\frac{3}{2}$};
\draw (0.03,1.0) -- (-0.03,1.00) node [left] {$1$};
\draw (0.03,.5) -- (-0.03,.5) node [left] {$\frac{1}{2}$};
\draw (0.03,-.5) -- (-0.03,-.5) node [left] {$-\frac{1}{2}$};
\draw (0.03,-1.0) -- (-0.03,-1.00) node [left] {$-1$};

% Plot
\draw[dashed] (1.0,0.0) -- (1.0,1.0);
\draw[fill=black!70, opacity=0.4] (0.0,-.5) -- (0.0,0.0) -- (.5,.5) -- (1.5,0.5)
-- (1.0,0.0) -- (.5,-.5) -- (0.0,-.5);
\draw (0.0,-1.0) -- (0.0,0.0) -- (1.0,1.0) -- (2.0,1.0) -- (1.0,0.0) --
(0.0,-1.0);
%\draw (0.85,0.65) node {};
%\draw (0.65,0.3) -- (1.15,0.8);
%\fill (1.15,0.8) circle[radius=1pt] node [right] {\small
%$F^{s_1}_{p_1,2}$} ;
%\fill (0.65,0.3) circle[radius=1pt] node [below] {\small
%$F^{s_0}_{p_0,2}$} ;
%\fill (0.925,0.575) circle[radius=1pt];
%\draw (0.15,-0.65) node {};
\end{tikzpicture}
\caption{Unconditional convergence 
%the Haar system in 
in Hardy-Sobolev spaces}\label{fig2}
\end{center}
\end{figure}

 In previous works stronger notions of convergence have been examined, such as   unconditional convergence for  the martingale difference series. This  is equivalent with the inequality
\Be \label{mult}\Big\|\sum_n b_n \bbD_n f\Big\|_{F^s_{p,q} } \lc \|b\|_{\ell^\infty(\bbN)} \|f\|_{F^s_{p,q}}.
\Ee 
It follows from the results in Triebel \cite{triebel-bases} that  \eqref{mult} holds if we add the condition
$1/q-1<s<1/q$ to the hypotheses in the theorem. For the case $q=2$ this corresponds to the shaded region in Figure \ref{fig2}.

It was shown in \cite{su}, \cite{sudet}  that the additional restriction on the $q$-parameter 
 is necessary for \eqref{mult} to hold. If we drop it  then Theorem 
\ref{expthm} and a summation by parts argument  imply  that \eqref{mult} holds with the larger norm $\|b\|_\infty +\|b\|_{BV}$.
It should be interesting to  establish  sharp  results involving sequence spaces that are intermediate between $\ell^\infty(\bbN)$ and $BV(\bbN)$. We remark that these problems are interesting only
for the $F^s_{p,q}$ spaces since  inequality \eqref{mult} with $B^s_{p,q}$ in place of $F^s_{p,q}$ holds in the full parameter range of Theorem \ref{expthm}, see   \cite{triebel-bases} for further discussion and historical comments.

%\Be\label{martdiff} \bbE_{N+1} f-\bbE_Nf= \sum_{\mu\in \bbZ} 2^N \inn{f}{ h_{N,\mu} }h_{N,\mu},\Eei.e.  $\bbE_{N+1} -\bbE_N $ is the orthogonal projection to the space generated by the Haar functions with Haar frequency $2^N$.Indeed it is easy to check from the definitions  thatfor each $N,\mu$ one has  $$[\bbE_{N+1} f-\bbE_Nf]\bbone_{I_{N,\mu}}=2^N \inn{f}{ h_{N,\mu} }h_{N,\mu}.$$The main tool in proving \cite[Thm.\ 1.1]{gsu-basic} is to prove the universal boundedness of the operators $\bbE_N$ in $F^s_{p,q}$. Note, that by earlier results of the last two authors \cite{su, sudet} it turned out that there is a ``solid'' parameter region in the 
%$(1/p,s)$-diagram, where the Haar system is a Schauder basis but not an unconditional basis in $F^s_{p,q}$. 

In \S\ref{proofsect} we give a proof of Theorem \ref{expthm} using characterizations of Triebel-Lizorkin spaces based on Daubechies wavelets. Relying on this, 
the proof is rather elementary due to the orthogonality and locality properties of the wavelet system. In addition, a ``wavelet analog'' of \cite[Thm.\ 1.2]{gsu-basic} is provided in Proposition \ref{ENPN} below. 
In \S\ref{Schauder} we apply the  methods to get an additional result needed to obtain the Schauder basis property of the Haar system.

\section{Proof of Theorem \ref{expthm}} \label{proofsect}
We will exclusively use a characterization of Triebel-Lizorkin spaces $F^s_{p,q}(\rr)$ and Besov spaces $B^s_{p,q}$ via compactly supported 
Daubechies wavelets \cite{Dau92}, \cite[Sect.\ 4]{Woj97}. Let $\psi_0$ and $\psi$ be the orthogonal scaling function and corresponding wavelet of Daubechies type such that 
$\psi_0, \psi$ being sufficiently smooth ($C^K$) and $\psi$ having sufficiently many vanishing moments ($L$). We denote 
$$
      \psi_{j,\nu}(x):=\frac{1}{\sqrt{2}}\psi(2^{j-1}x-\nu)\quad,\quad j\in \bbN, \nu \in \bbZ\,,
$$
and $\psi_{0,\nu}(x):=\psi_0(x-\nu)$ for $\nu \in \bbZ$. Let $0<p<\infty$, $0<q\leq \infty$ and $s\in \bbR$. If $K$ and $L$ are large enough (depending on $p,q$ and $s$)
then we have the equivalent characterization (usual modification in case $q=\infty$), 
\begin{eqnarray}\label{wave_charF}
      \|f\|_{F^s_{p,q}} &\asymp& \Big\|\Big(\sum\limits_{j=0}^{\infty} \Big|2^{js} \sum\limits_{\nu \in \bbZ} \lambda_{j,\nu}(f)\bbone_{j,\nu}\Big|^q\Big)^{1/q}\Big\|_p\,,\\
    \label{wave_charB}\|f\|_{B^s_{p,q}} &\asymp& \Big(\sum\limits_{j=0}^{\infty} \Big\|2^{js} \sum\limits_{\nu \in \bbZ} \lambda_{j,\nu}(f)\bbone_{j,\nu}\Big\|_p^q\Big)^{1/q}\,,
  \end{eqnarray}
where $\lambda_{j,\nu}(f):=2^j\langle f, \psi_{j,\nu} \rangle$ and $\bbone_{j,\nu}$ denotes the characteristic function of the interval $I_{j,\nu}:=[2^{-j}\nu,2^{-j}(\nu+1)]$.
See Triebel \cite[Thm.\ 1.64]{Tr06} and the references therein.
A corresponding characterization also holds true for Besov spaces $B^s_{p,q}$. Since we also 
deal with distributions which are not locally integrable, the inner product $\langle f,\psi_{j,\nu} \rangle$ has to be interpreted in the usual way. Clearly, 
$f$ can be decomposed into wavelet building blocks, i.e. 
\begin{equation}\label{decomp}
    f = \sum\limits_{j\in \bbZ} f_j \, \text{ with } \, 
    f_j= \begin{cases} \sum\limits_{\nu\in\bbZ} \lambda_{j,\nu}(f) \psi_{j,\nu}\,,&\text{ if } j\ge 0,
    \\
    0 &\text{ if } j< 0.
    \end{cases}
\end{equation}
Let us denote the $N$th partial sum of this representation by 
\begin{equation}\label{def:PN}
    P_Nf = \sum\limits_{j \leq N} f_j\quad,\quad N\in \mathbb{N}\,.
\end{equation}

%where we simply put $f \equiv 0$ if $j<0$\,. 
Note, that the functions $f_j$ and $P_Nf$ represent $K$ times 
continuously differentiable functions due to the regularity assumption on the wavelet. 

In the sequel we will prove the following proposition. 
\begin{proposition}\label{ENPN} Let $1/2<p\leq \infty$, $0<r\leq\infty$, and $1/p-1<s<\min\{1/p,1\}$. 
Let $\{\psi_{j,\nu}\}_{j,\nu}$ represent a Daubechies wavelet system such that \eqref{wave_charB} holds for all $0<q\leq \infty$ and 
let $P_N$ be given by \eqref{def:PN}. Then there is a constant $C:=C(p,r,s)>0$ such that for all $f\in B^s_{p,\infty}$\\
\begin{equation}
    \sup\limits_{N\in \bbN}\|\bbE_Nf-P_N f\|_{B^s_{p,r}} \leq C\|f\|_{B^s_{p,\infty}}\,.
\end{equation}
\end{proposition}
Note, that for a fixed wavelet system satisfying \eqref{wave_charF}  we clearly have 
\begin{equation}\label{PNbound}
  \sup\limits_{N\in \bbN}\|P_N f\|_{F^s_{p,q}} \leq C\|f\|_{F^s_{p,q}}\,.
\end{equation}
If this wavelet system in addition satisfies \eqref{wave_charB} for all $0<q\leq \infty$ then 
Proposition \ref{ENPN} together with \eqref{PNbound} implies Theorem \ref{expthm}. 

 \begin{proof}[\bf A. Proof of Proposition \ref{ENPN} in the case $\mathbf{1/2 \boldsymbol < p\boldsymbol\le 1}$.] 
 {Let $1/p-1<s<1$}. 
 %\subsection{Proof   in the case  $1/2<p\leq 1$} \label{qBcase} 
 %Then there is a universal constant $C>0$ such that  \begin{proposition}\label{qBanach} Let $1/2<p\leq 1$, $0<q\leq \infty$ and $1/p-1<s<1$. Then there is a universal constant $C>0$ such that for all $N\in \bbN$ and all $f\in F^s_{p,q}(\rr)$$$    \|\bbE_N(f)\|_{F^s_{p,q}} \leq C\|f\|_{F^s_{p,q}}\,.$$\end{proposition}
Using the decomposition \eqref{decomp} 
we can write with $\theta:=\min\{1,p\} = p$
\begin{eqnarray}\label{eq1}
    \nonumber\|\bbE_N f-P_Nf\|_{B^s_{p,r}} &\asymp& \Big(\sum\limits_{j=0}^{\infty}\Big\|2^{js}
    \sum\limits_{\eta}2^j\langle \bbE_N\!f-P_Nf, \psi_{j,\eta}\rangle \bbone_{j,\eta}
    \Big\|_p^r\Big)^{1/r}\\
    \label{eq1_a}&\lesssim& \Big(\sum\limits_{j=0}^{\infty}\Big\|2^{js}\sum\limits_{\eta}2^j\langle \bbE_N\!(P_Nf)-P_Nf, \psi_{j,\eta}\rangle \bbone_{j,\eta}
    \Big\|_p^r\Big)^{1/r}\\
    \label{eq1_a1}&&\hspace{0.5cm}+\Big(\sum\limits_{j=0}^{\infty}\Big\|2^{js}
    \sum\limits_{\eta}2^j\langle \bbE_N(f-P_Nf), \psi_{j,\eta}\rangle \bbone_{j,\eta}\Big\|_p^r\Big)^{1/r}\,.
\end{eqnarray}
We split the proof into several steps according to the cases we have to distinguish in the estimation of the 
quantities in \eqref{eq1} and \eqref{eq1_a1}.
 \smallskip
 
\indent {\em Step A1.} We deal with \eqref{eq1_a1} and use that $f-P_Nf = \sum\limits_{j+\ell>N} f_{j+\ell}$. Clearly,\\
\begin{equation}\label{thetatriang}
    \eqref{eq1_a1} \lesssim \Big(\sum\limits_{j= 0}^{\infty}\Big(\sum\limits_{j+\ell\geq N}\Big\|2^{js}\sum\limits_{\eta}2^j\langle \bbE_N f_{j+\ell}, \psi_{j,\eta}\rangle \bbone_{j,\eta}
    \Big\|_p^\theta\Big)^{r/\theta}\Big]^{1/r}
\end{equation}
We continue estimating $\|2^{js}\sum_{\eta}2^j\langle \bbE_N f_{j+\ell}, \psi_{j,\eta}\rangle \bbone_{j,\eta}
    \|_p$. Note first that due to $p\leq 1$
\begin{equation}\label{wav_red}
\begin{split}
&\Big\|2^{js}\sum_{\eta}2^j\langle \bbE_N f_{j+\ell}, \psi_{j,\eta}\rangle \bbone_{j,\eta}\Big\|_p\\
&\hspace{0.5cm}\leq \Big(\sum\limits_{\nu \in \bbZ} |\lambda_{j+\ell,\nu}(f)|^p\Big\|2^{js}\sum_{\eta}2^j\langle \bbE_N \psi_{j+\ell,\nu}, \psi_{j,\eta}\rangle \bbone_{j,\eta}\Big\|_p^p\Big)^{1/p}\,.
\end{split}
\end{equation}
So it remains to deal with $\|2^{js}\sum_{\eta}2^j\langle \bbE_N \psi_{j+\ell,\nu}, \psi_{j,\eta}\rangle \bbone_{j,\eta}\|_p$. Note, that due to $j+\ell> N$
the function $\bbE_N\psi_{j+\ell,\nu}$ is a step function consisting of $O(1)$ non-vanishing steps. These steps have length $2^{-N}$ and 
magnitude bounded by $O(2^{N-(j+\ell)})$.

\indent {\em Case A1.1} Assume $j\geq N$.\\ Due to the cancellation of $\psi_{j,\eta}$ and $j\geq N$ we have that the function
$\sum_{\eta}2^j\langle \bbE_N \psi_{j+\ell,\nu}, \psi_{j,\eta}\rangle \bbone_{j,\eta}$ 
is supported on a union of intervals of total measure $O(2^{-j})$ and bounded from above by $O(2^{N-(j+\ell)})$. This gives
\begin{equation}\label{eq11}
    \Big\|2^{js}\sum_{\eta}2^j\langle \bbE_N \psi_{j+\ell,\nu}, \psi_{j,\eta}\rangle \bbone_{j,\eta}\Big\|_p \lesssim 2^{js}2^{-j/p}2^{N-j-\ell}\,.
\end{equation}
\indent {\em Case A1.2.} Assume $j \leq N$.\\
Clearly, we have $\ell > 0$ since $j+\ell> N$. 
%By straightforward size estimates
Now 
$\sum_{\eta}2^j\langle \bbE_N \psi_{j+\ell,\nu}, \psi_{j,\eta}\rangle \bbone_{j,\eta}$ 
is supported on an interval of size $O(2^{-j})$. 
%Since $j \leq  N$ we can not make use of the cancellation properties of $\varphi_j$. 
%We still have that $|\bbE_N(\psi_{j+\ell,\nu})\ast \varphi_j(x)|$
As $\bbE_N\psi_{j+\ell,\nu}$ consists of $O(1)$ steps of length $2^{-N}$ each and $N\geq j$ 
we get
by straightforward size estimates
$2^j\langle\bbE_N\psi_{j+\ell,\nu},\psi_{j,\eta}\rangle=O(2^{-\ell})$. Hence
%the convolution produces an additional factor $O(2^{j-N})$. Hence,
\begin{equation}\label{eq12}
   \Big\|2^{js}\sum_{\eta}2^j\langle \bbE_N \psi_{j+\ell,\nu}, \psi_{j,\eta}\rangle \bbone_{j,\eta}\Big\|_p  \lesssim 2^{js}2^{-j/p}2^{-\ell}\,.
\end{equation}

\smallskip
\indent{\em Step A2.} We  consider \eqref{eq1_a} and observe first 
\begin{equation}\label{9_1}
    \eqref{eq1_a} \lesssim \Big(\sum\limits_{j= 0}^{\infty}
    \Big(\sum\limits_{j+\ell\leq N}\Big\|2^{js}\sum\limits_{\eta}2^j\langle \bbE_N f_{j+\ell}-f_{j+\ell}, \psi_{j,\eta}\rangle\bbone_{j,\eta}
    \Big\|_p^\theta\Big)^{r/\theta}\Big)^{1/r}\,.
\end{equation}
Analogously to \eqref{wav_red} the matter reduces to estimate the $L_p$ \\(quasi-)norm of the functions
\begin{equation}\label{reduce}
2^{js}\sum\limits_{\eta}2^j\langle \bbE_N \psi_{j+\ell,\nu}-\psi_{j+\ell,\nu}, \psi_{j,\eta}\rangle\bbone_{j,\eta}
\end{equation}
for the different cases resulting from $j+\ell\leq N$\,.\\
\indent {\em Case A2.1.} %Assume $\ell \geq 0$ and 
We first deal with the case $j \leq N$. 
%and estimate as follows
% %\begin{equation}\label{eq31}
% \begin{subequations}
%  \begin{align}
%     &\Big\|\Big(\sum\limits_{j}|2^{js}\vpj*\bbE_N\!f_{j+\ell}|^q\Big)^{1/q}\Big\|^{\theta}_p\notag\\
%     &\lesssim 
%     \Big\|\Big(\sum\limits_{j}|2^{js} \vpj*[\bbE_N\!f_{j+\ell}-f_{j+\ell}]|^q\Big)^{1/q}\Big\|^{\theta}_p
%     \label{eq31}\\
%     &\hspace{0.5cm}+ \Big\|\Big(\sum\limits_{j}|2^{js}\vpj*f_{j+\ell}|^q\Big)^{1/q}\Big\|^{\theta}_p\,.
%     \label{eq31sec}
%  \end{align}
% \end{subequations}
% Similar as in \eqref{thetatriang} we estimate the  term in  \eqref{eq31} via 
% \begin{equation}\label{eq31a}
%   \begin{split}
%    &\Big\|\Big(\sum\limits_{j}^{\infty}|2^{js}\vpj*[\bbE_N\!f_{j+\ell}-f_{j+\ell}]|^q\Big)^{1/q}\Big\|^{\theta}_p\\
%    &\hspace{0.5cm}\leq \sum\limits_{j}\|2^{js}\vpj*[\bbE_N\!f_{j+\ell}-f_{j+\ell}]\|_p^{\theta}\,.
%   \end{split} 
% \end{equation}
% Again, analogously to \eqref{wav_red} we have 
% \begin{equation}\label{eq31b}
%   \begin{split}
%   &\Big\|2^{js}\sum\limits_{\eta}2^{j}\langle \bbE_N\!f_{j+\ell}-f_{j+\ell},\psi_{j,\eta}\rangle \bbone_{j,\eta}\Big\|_p\\
%     &\hspace{0.5cm}\lesssim \Big(\sum\limits_{\nu \in \bbZ} |\lambda_{j+\ell,\nu}(f)|^p
%     \Big\|2^{js}\sum\limits_{\eta}2^{j}\langle \bbE_N\psi_{j+\ell,\nu}-\psi_{j+\ell,\nu},\psi_{j,\eta}\rangle \bbone_{j,\eta}\Big\|_p^p\Big)^{1/p}\,.
%   \end{split}
% \end{equation}
Using the mean value theorem 
%of calculus
 together with \eqref{condexp} we see for all $x\in \rr$ that 
$$
%   |\vpj*[\bbE_N\psi_{j+\ell,\nu}-\psi_{j+\ell,\nu}](x)| \leq 2^{j+\ell-N}\,.
   |\bbE_N\psi_{j+\ell,\nu}(x)-\psi_{j+\ell,\nu}(x)| \leq 2^{j+\ell-N}\,.
$$
Due to $j+\ell\leq N$, its support has length $O(2^{-(j+\ell)})$ around $\nu 2^{-(j+\ell)}$.
We continue distinguishing the  cases $\ell\ge 0$ and $\ell<0$. \\
\indent {\em Case A2.1.1.} Let $\ell \geq 0$. Since $j+\ell\geq j$ the inner product with $2^j\psi_{j,\eta}$ 
gives an additional factor $2^{-\ell}$. In addition, the support of 
\eqref{reduce} is contained in 
an interval of size $O(2^{-j})$. Hence, 
we get 
\begin{equation}\label{eq311}
  \Big\|2^{js}\sum\limits_{\eta}2^{j}\langle \bbE_N\psi_{j+\ell,\nu}-\psi_{j+\ell,\nu},\psi_{j,\eta}\rangle \bbone_{j,\eta}\Big\|_p\lesssim 
  2^{js}2^{j+\ell-N}2^{-\ell}2^{-j/p}\,.
\end{equation}
\indent {\em Case A2.1.2.} Assume $\ell \leq 0$. This time the inner product with $2^j\psi_{j,\eta}$ does not give an extra factor and the support 
has length $2^{-(j+\ell)}$. Thus, we have in this case
\begin{equation}\label{eq312}
    \Big\|2^{js}\sum\limits_{\eta}2^{j}\langle \bbE_N\psi_{j+\ell,\nu}-\psi_{j+\ell,\nu},\psi_{j,\eta}\rangle \bbone_{j,\eta}\Big\|_p\lesssim 2^{js}2^{j+\ell-N}2^{-(j+\ell)/p}\,.
\end{equation}

\indent {\em Case A2.2.} Assume $j > N \geq j+\ell$ which implies $\ell <0$. 
Due to the orthogonality of the wavelets ($\ell<0$)
we can estimate
%the following identity
\begin{equation}\label{eq32}
   \begin{split}
	&\Big\|2^{js}\sum\limits_{\eta}2^{j}\langle \bbE_N\psi_{j+\ell,\nu}-\psi_{j+\ell,\nu},\psi_{j,\eta}\rangle \bbone_{j,\eta}\Big\|_p\\
	&~~~~~\lesssim 2^{js}\Big(\sum\limits_{\mu \in \bbZ} \int\limits_{|x-2^{-N}\mu|\lesssim 2^{-j}}\Big|
	\sum\limits_{\eta}2^j\langle \bbE_N\psi_{j+\ell,\nu},\psi_{j,\eta}\rangle\bbone_{j,\eta}(x)\Big|^p\,dx\Big)^{1/p}\\
	&~~~~~\lesssim 2^{js}\Big(\sum\limits_{\mu \in \bbZ} \int\limits_{|x-2^{-N}\mu|\lesssim 2^{-j}}\Big|
	\sum\limits_{\eta}2^j\langle \bbE_N\psi_{j+\ell,\nu}-\psi_{j+\ell,\nu},\psi_{j,\eta}\rangle\bbone_{j,\eta}(x)\Big|^p\,dx\Big)^{1/p}\\
	&~~~~~\lesssim 2^{js}2^{j+\ell-N}2^{[N-(j+\ell)-j]/p}\,,
   \end{split}
\end{equation}
where we took into account that the $\mu$-sum consists of $O(2^{N-(j+\ell)})$ summands. 
\smallskip
%{\bf Hieraus ergibt sich: $s<1$}\\

\indent {\em Step A3. Estimation of \eqref{eq1_a1}.} Plugging \eqref{wav_red} and \eqref{eq11} into the right hand side of \eqref{thetatriang} yields 
\begin{equation}\label{HoelNik}
  \begin{split}
   &\Big[\sum\limits_{j= N}^{\infty}\Big(\sum\limits_{j+\ell\geq N}\Big\|2^{js}\sum\limits_{\eta}2^j\langle \bbE_N f_{j+\ell}, \psi_{j,\eta}\rangle \bbone_{j,\eta}
    \Big\|_p^\theta\Big)^{r/\theta}\Big]^{1/r}\\
    &\lesssim
    A_N\sup\limits_{j,\ell} 
    \Big(\sum\limits_{\nu \in \bbZ} |2^{(j+\ell)s}\lambda_{j+\ell,\nu}(f)|^p2^{-(j+\ell)}\Big)^{1/p}\\
    &\lesssim A_N\|f\|_{B^s_{p,\infty}}
  \end{split}  
\end{equation}
with
$$
A_N^r=
\sum\limits_{j\geq N}2^{(N-j)r}\Big(\sum\limits_{\ell\geq N-j}2^{\theta\ell(1/p-1-s)}\Big)^{r/\theta}\lesssim 1 
$$
%$1/p-1<s<1/p$ then the sums are uniformly bounded (in $N$).
by the assumption $1/p>s>1/p-1$. \\
\indent Plugging \eqref{wav_red} and \eqref{eq12} into into the right hand side of \eqref{thetatriang} leads to a similar estimate as above, only the sums over $j$ and $\ell$ change 
to 
\[
  \widetilde A_N^r= \sum\limits_{j\leq N}\Big(\sum\limits_{\ell\geq N-j}2^{\theta\ell(1/p-1-s)}\Big)^{r/\theta}
\]
which is uniformly bounded in $N$ if $s>1/p-1$\,.

\smallskip
{\em Step A4. Estimation of \eqref{eq1_a}.} Combining \eqref{9_1}, \eqref{eq311} and \eqref{eq312} we find
\begin{equation}\nonumber
  \begin{split}
   &\Big[\sum\limits_{j= 0}^{N}\Big(\sum\limits_{j+\ell\leq N}\Big\|
   2^{js}\sum\limits_{\eta}2^j\langle \bbE_N f_{j+\ell}-f_{j+\ell}, \psi_{j,\eta}\rangle \bbone_{j,\eta}\Big\|_p^\theta\Big)^{r/\theta}\Big]^{1/r}\\
    &\lesssim \Big[\Big(\sum\limits_{j\leq N}2^{(j-N)r}
    \Big(\sum\limits_{\ell=-\infty}^{N-j}2^{\theta\ell(1/p-s)}\Big)^{r/\theta}\Big]^{1/r}
    \|f\|_{B^s_{p,\infty}}\,.
  \end{split}  
\end{equation}
The sums are finite and uniformly bounded if $1/p-1<s<1/p$.\\
Finally, we combine \eqref{thetatriang}, \eqref{wav_red} and \eqref{eq32} to obtain
\begin{equation}\label{casej>N}
 \begin{split}
   &\Big(\sum\limits_{j= N}^{\infty}\Big(\sum\limits_{j+\ell\leq N}\Big\|
   2^{js}\sum\limits_{\eta}2^j\langle \bbE_N f_{j+\ell}-f_{j+\ell}, \psi_{j,\eta}\rangle \bbone_{j,\eta}\Big\|_p^\theta\Big)^{r/\theta}\Big]^{1/r}\\
    &\lesssim \Big[\Big(\sum\limits_{j\geq N}2^{(j-N)\theta}2^{(N-j)\theta/p}\Big(\sum\limits_{\ell=-\infty}^{N-j}2^{r\ell(1-s)}\Big)^{r/\theta}\Big]^{1/r}
    \|f\|_{B^s_{p,\infty}}\,,
 \end{split}
\end{equation}
which is uniformly bounded if $s<1$. This concludes the proof in the case $p\le 1$.
\end{proof}
\begin{proof}[\bf B. Proof in the case $\mathbf{1 \leq \boldsymbol p\boldsymbol{\leq \infty}}$.]
%of Theorem \ref{expthm}, conclusion}
%\vspace{-0.3cm}
%\begin{proof} 
%It remains to prove the theorem in 
%Taking Proposition \ref{qBanach} into account it remains to deal with
 %the cases 
%$1 < p<\infty$, $0<q\leq \infty$. \\
We follow the proof in the case $p\le 1$ 
%Proposition \ref{qBanach} 
%and follow that proof 
until
\eqref{thetatriang} and \eqref{9_1}, respectively. Note, that we may use $\theta = 1$ now. 
Then we have to proceed differently. 

\indent {\em Case B1.1}  Assume $N < j, j+\ell$. 
Taking \eqref{thetatriang} into account we replace \eqref{wav_red} by 
\begin{equation}\label{p>1}
 \begin{split}
    &\Big\|2^{js}\sum\limits_{\eta}2^j\langle \bbE_N\!f_{j+\ell},\psi_{j,\eta}\rangle \bbone_{j,\eta}\Big\|^p_p \\
    &\hspace{0.5cm}\leq \int\Big[\sum\limits_{\nu \in \bbZ}|2^{js}\lambda_{j+\ell,\nu}(f)|\cdot \Big|\sum\limits_{\eta}2^j\langle 
    \bbE_N\!\psi_{j+\ell,\nu},\psi_{j,\eta}\rangle \bbone_{j,\eta}(x)\Big|\Big]^p\,dx\\
    &\hspace{0.5cm}\lesssim \sum\limits_{\nu \in \bbZ} |2^{js}\lambda_{j+\ell,\nu}(f)|^p2^{-j}\, 2^{(N-j-\ell)p}.
 \end{split}
\end{equation}
Indeed, since $\bbE_N\psi_{j+\ell,\nu}= 0$ if $\supp\,\psi_{j+\ell,\nu} \subset I_{N,\mu}$ 
the sum on the right-hand side of \eqref{p>1} is lacunary and the 
%appearing
 functions $\sum_{\eta}2^j\langle 
    \bbE_N\!\psi_{j+\ell,\nu},\psi_{j,\eta}\rangle \bbone_{j,\eta}$ have essentially disjoint support 
    (for different $\nu$). Hence, we get 
\begin{equation}\label{20_1}
 \begin{split}
   &\Big\|2^{js}\sum\limits_{\eta}2^j\langle \bbE_N\!f_{j+\ell},\psi_{j,\eta}\rangle \bbone_{j,\eta}\Big\|_p \\
   &\hspace{0.5cm}\lesssim 2^{-\ell s}2^{N-j-\ell}2^{\ell/p}\Big(\sum\limits_{\nu \in \bbZ}|2^{(j+\ell)s}\lambda_{j+\ell,\nu}(f)|^p2^{-(j+\ell)}\Big)^{1/p}\\
   &\hspace{0.5cm}\lesssim2^{-\ell s}2^{N-j-\ell}2^{\ell/p}\|f\|_{B^s_{p,\infty}}\,.
   \end{split}
\end{equation}
For $1/p-1<s<1/p$ the sum over the respective range of $j$ and $\ell$ is uniformly bounded.
\\
\indent {\em Case B1.2.} We  now deal with $j+\ell> N \geq j$. Due to the orthogonality of the wavelet system and $\ell>0$ we obtain
\begin{equation}\label{20_2}
    \Big\|2^{js}\sum\limits_{\eta}2^j\langle \bbE_N\!f_{j+\ell},\psi_{j,\eta}\rangle \bbone_{j,\eta}\Big\|_p = 
    \Big\|2^{js}\sum\limits_{\eta}2^j\langle \bbE_N\!f_{j+\ell}-f_{j+\ell},\psi_{j,\eta}\rangle \bbone_{j,\eta}\Big\|_p\,.
\end{equation}
We continue exploiting the cancellation property
\begin{equation}\label{canc}
    %0 = \bbE_N(\bbE_N(f)-f)(x) = \int_{I_{N,\mu(x)}} \bbE_N(f)(y)-f(y)\,dy\,.
    \bbE_N(f-\bbE_N\,f)=0\,
\end{equation}
to estimate the right-hand side of \eqref{20_2}. We obtain the following identities
%obtain the pointwise estimate
%\begin{equation}\label{mu-conv}
% \begin{split}   &\Big|2^{js}\int\varphi_j(x-y)(\bbE_N\!f_{j+\ell}(y)-f_{j+\ell}(y))\,dy\Big|\\   &\hspace{0.5cm}=\Big|2^{js}\sum\limits_{\mu:|2^{-N}\mu-x|\lesssim 2^{-j}}\int_{I_{N,\mu}}\varphi_j(x-y)(\bbE_N\!f_{j+\ell}(y)-f_{j+\ell}(y))\,dy\Big|\\   &\hspace{0.5cm}=\Big|2^{js}\sum\limits_{\mu:|2^{-N}\mu-x|\lesssim 2^{-j}}\int_{I_{N,\mu}}(\varphi_j(x-y)-\varphi_j(x))\times\\   &\hspace{4.5cm}(\bbE_N\!f_{j+\ell}(y)-f_{j+\ell}(y))\,dy\Big|\, \end{split}\end{equation}
\begin{align}\label{mu-conv}
% \begin{split}
   &\Big|2^{js}\sum_{\eta}\bbone_{j,\eta}(x)2^j\int \psi_{j,\eta}(y)(\bbE_N\!f_{j+\ell}(y)-f_{j+\ell}(y))\,dy\Big|\\
   \notag
   &=\Big|2^{js}\sum_{\eta}\bbone_{j,\eta}(x)\sum\limits_{\mu:|2^{-N}\mu-x|\lesssim 2^{-j}}2^j\int_{I_{N,\mu}}\psi_{j,\eta}(y)(\bbE_N\!f_{j+\ell}(y)-f_{j+\ell}(y))\,dy\Big|\\
   \notag
   &=\Big|2^{js}\sum_{\eta}\bbone_{j,\eta}(x)\sum\limits_{\mu:|2^{-N}\mu-x|\lesssim 2^{-j}}\\
   &\hspace{0.6cm}2^j\int_{I_{N,\mu}}(\psi_{j,\eta}(y)-\psi_{j,\eta}(2^{-N}\mu))(\bbE_N\!f_{j+\ell}(y)-f_{j+\ell}(y))\,dy\Big|\, \notag\,.
 \end{align}
Let $\eta \in \bbZ$ such that $\bbone_{j,\eta}(x) = 1$. We continue estimating \eqref{mu-conv} by
%estimating
\begin{equation}\nonumber
 \begin{split}
   &\hspace{0.5cm} 2^{js}\sum\limits_{\mu:|2^{-N}\mu-x|\lesssim 2^{-j}}2^j\int_{I_{N,\mu}}|
   (\psi_{j,\eta}(y)-\psi_{j,\eta}(2^{-N}\mu))\cdot \bbE_N\!f_{j+\ell}(y)|\,dy\\
   &\hspace{0.9cm}+\Big|2^{js}\sum\limits_{\mu:|2^{-N}\mu-x|\lesssim 2^{-j}}2^j\int_{I_{N,\mu}}
   (\psi_{j,\eta}(y)-\psi_{j,\eta}(2^{-N}\mu))\cdot f^{\mu,1}_{j+\ell}(y)\,dy\Big|\\
   &\hspace{0.9cm}+\Big|2^{js}\sum\limits_{\mu:|2^{-N}\mu-x|\lesssim 2^{-j}}2^j\int_{I_{N,\mu}}
   (\psi_{j,\eta}(y)-\psi_{j,\eta}(2^{-N}\mu))\cdot f^{\mu,2}_{j+\ell}(y)\,dy\Big|\\
   &\hspace{0.5cm}=:F_0(x)+F_1(x)+F_2(x)\,,
   \end{split}
\end{equation}
where 
\begin{equation}
 \begin{split}
    \nonumber f_{j+\ell}^{\mu}&:=\sum\limits_{\nu:\supp\,\psi_{j+\ell,\nu} \cap I_{N,\mu} \neq \emptyset}\lambda_{j+\ell,\nu}(f)\psi_{j+\ell,\nu}\,,\\
    \nonumber f_{j+\ell}^{\mu,1}&:=\sum\limits_{\nu:\supp\,\psi_{j+\ell,\nu} \subset I_{N,\mu}}\lambda_{j+\ell,\nu}(f)\psi_{j+\ell,\nu}\,,
    \\f_{j+\ell}^{\mu,2}&:=f_{j+\ell}^{\mu}-f_{j+\ell}^{\mu,1}\,.
    \nonumber
 \end{split}
\end{equation}
%$f_{j+\ell}^{\mu,2}:=f_{j+\ell}^{\mu}-f_{j+\ell}^{\mu,1}$,
% and $g_{j+\ell}^{1}:=\sum_{\mu} f_{j+\ell}^{\mu,1}$. 
Note, that the function $F_1$ vanishes since $\ell> 0$ (use orthogonality) and $j+\ell>0$ (use vanishing moments).\\
$F_0(x)$ can be  estimated by
% from above by
$$
    2^{js}\sum\limits_{\mu:|2^{-N}\mu-x|\lesssim 2^{-j}}2^{2j-2N}\sup\limits_{y\in I_{N,\mu}} \sum\limits_{\nu:\supp\,\psi_{j+\ell,\nu} \cap I_{N,\mu} \neq \emptyset}|
    \lambda_{j+\ell,\nu}(f) \bbE_N(\psi_{j+\ell,\nu})(y)|\,.
$$
Here $\bbE_N\psi_{j+\ell,\nu}$ is mostly vanishing, namely when $\supp\,\psi_{j+\ell,\nu} \subset I_{N,\mu}$. If it does not vanish then the boundary of 
$I_{N,\mu}$ intersects $\supp\,\psi_{j+\ell,\nu}$ 
%is contained in its support
 and $|\bbE_N\psi_{j+\ell,\nu}| \lesssim 2^{N-(j+\ell)}$. This happens only for a bounded number of $\nu$'s (independently  of $j,\ell$). Thus for a fixed $y$  only a bounded number of coefficients contribute. Hence, we have
\begin{equation}\label{eq27a}
% \begin{split}
 F_0(x)\lesssim 2^{js}2^{2j-2N}2^{N-(j+\ell)}
  %\times\\  &\times
  \sum\limits_{\mu:|2^{-N}\mu-x|\lesssim 2^{-j}} \sup\limits_{\nu:\supp\,\psi_{j+\ell,\nu} \cap\partial  I_{N,\mu} \neq \emptyset}|\lambda_{j+\ell,\nu}(f)|\,.
 %\end{split}
\end{equation}
Taking the $L_p$-norm and using H\"older's inequality with $1/p+1/p'=1$ yields
\begin{multline}\label{F0}
% \begin{split}  &
\|F_0\|_p\lesssim 2^{-\ell s}2^{2j-2N}2^{N-(j+\ell)}2^{(N-j)/p'}2^{\ell/p}\times\\
  %&\hspace{1.5cm}
  \Big(\sum\limits_{\nu}|2^{(j+\ell)s}\lambda_{j+\ell,\nu}(f)|^p2^{-(j+\ell)}\Big)^{1/p}\,,
% \end{split}
\end{multline}
where again $\|f\|_{B^s_{p,\infty}}$ dominates the sum on the right-hand side, see \eqref{wave_charB}.
Finally, we deal with $F_2(x)$. Since to $f^{\mu,2}_{j+\ell}$ only a uniformly bounded number of coefficients $\lambda_{j+\ell,\nu}$ contribute 
to the sum and the integrals are  
taken over an interval of length $O(2^{-(j+\ell)})$ we obtain, similar as above, by H\"older's inequality 
\begin{equation}\label{F_3}
    \|F_2\|_p \lesssim 2^{-\ell s}2^{-\ell+j-N}2^{(N-j)/p'}2^{\ell/p}\Big(\sum\limits_{\nu}|2^{(j+\ell)s}\lambda_{j+\ell,\nu}(f)|^p2^{-(j+\ell)}\Big)^{1/p}\,.
\end{equation}
%\footnote{GG: should be $2^{-\ell} 2^{j-N}$ if it comes from MVT in $\phi$ and $2^{-j-\ell} $ due to the support of $\psi_{j+\ell,\nu}$}
Putting the estimates from \eqref{20_2} to \eqref{F_3} together we observe that the sum over the respective range of $j$ and $\ell$ (see \eqref{eq1_a1}) is uniformly bounded with 
respect to $N$ if $s>1/p-1$. 

\indent {\em Case B2.1.} Here we deal with $j+\ell,j\leq N$. Starting from \eqref{9_1} (with $\theta = 1$)
we continue similarly as after
\eqref{canc} and obtain the pointwise identity \eqref{mu-conv}. Note, that we already start with 
$\bbE_N f_{j+\ell}-f_{j+\ell}$, so we do have to use the orthogonality argument \eqref{20_2}, which does indeed not apply 
here since $\ell = 0$ is admitted. 

Since $j+\ell\leq N$ there is only a bounded number of coefficients 
$\lambda_{j+\ell,\nu}(f)$ contributing to $f_{j+\ell}$ on $I_{N,\mu}$. Using the mean value theorem in both factors of the 
integral in \eqref{mu-conv} we obtain 
\begin{multline*}
%\nonumber
 % \begin{split}
   \Big|2^{js}\sum\limits_{\eta}2^j \langle \bbE_N(f_{j+\ell})-f_{j+\ell},\psi_{j,\eta}\rangle \bbone_{j,\eta}\Big| \lesssim 2^{js}2^{2j-2N}2^{j+\ell-N}\times\\
   \sum\limits_{\mu:|2^{-N}\mu-x|\lesssim 2^{-j}} \sup\limits_{|\nu 2^{-(j+\ell)}-2^{-N}\mu|\lesssim 1} |\lambda_{j+\ell,\nu}(f)|\,,
%  \end{split}
\end{multline*}
which yields
\begin{equation}\nonumber
 \begin{split}
  &\Big|2^{js}\sum\limits_{\eta}2^j \langle \bbE_N(f_{j+\ell})-f_{j+\ell},\psi_{j,\eta}\rangle \bbone_{j,\eta}\Big\|_p \\
  &\hspace{0.5cm} \lesssim 2^{-\ell s}2^{j+\ell-N}2^{2j-2N}2^{(N-j)/p'}2^{\ell/p}\Big(\sum\limits_{\nu \in \bbZ}|2^{(j+\ell)s}\lambda_{j+\ell,\nu}(f)|^p2^{-(j+\ell)}\Big)^{1/p}\,.
 \end{split}
\end{equation}
The sum over the respective $j$ and $\ell$ is uniformly bounded in $N$ whenever $-1<s<1+1/p$.

{\em Case B2.2.} Finally $j+\ell\leq N < j$. Using again the orthogonality relation of the wavelets
we may estimate as follows (similar to \eqref{eq32})
\begin{equation}\label{final}
 \begin{split}
 &\Big\|2^{js}\sum\limits_{\eta}2^j\langle \bbE_N\!f_{j+\ell}-f_{j+\ell},\psi_{j,\eta}\rangle \bbone_{j,\eta}\Big\|_p\\
 &\lesssim 2^{js}\Big(\sum\limits_{\mu\in \bbZ} \int\limits_{|x-2^{-N}\mu|\lesssim 2^{-j}}
    \Big|\sum\limits_{\eta}2^j\langle \bbE_N\!f_{j+\ell}-f_{j+\ell},\psi_{j,\eta}\rangle \bbone_{j,\eta}(x)\Big|^p\, dx\Big)^{1/p}\,,
 \end{split}
\end{equation}
which is bounded by (see \eqref{eq32})
\begin{equation}\label{first}
   2^{-\ell s} 2^{j+\ell-N}2^{(N-j)/p}\Big(\sum\limits_{\nu \in \bbZ}|2^{(j+\ell)s}\lambda_{j+\ell,\nu}(f)|^p2^{-(j+\ell)}\Big)^{1/p}\,.
\end{equation}
Altogether we encounter the condition $1/p-1<s<1/p$ for any $0<r\leq \infty$ 
for the uniform boundedness of $\bbE_N:B^s_{p,\infty}\to B^s_{p,r}$ in case $1\leq p \leq \infty$.
\end{proof}

\section{On the Schauder basis property for the Haar system.
} \label{Schauder}
Let $\{h_{N,\mu}:\mu\in \bbZ\}$ be the set of Haar functions with Haar frequency  $2^{-N}$ and define 
for  {$N\in\bbN_0$}  and sequences  $a\in \ell^\infty(\bbZ)$,
\Be \label{TNdef}T_{N}[f,a]
=
\sum_{\mu\in \bbZ} a_\mu 2^{N}\inn{f}{h_{N,\mu}} h_{N,\mu}.
\Ee
In particular for the choice of $a=(1,1,1,\dots)$ one  recovers the operator $\bbE_{N+1}-\bbE_N$. %We have the following theorems. 
It was shown in \cite{gsu-basic} that Theorem \ref{expthm} together with 
%the inequality
\begin{multline}\label{T_N}
    \sup\limits_{N\in \bbN} \sup\limits_{\|a\|_{\infty}\leq 1} \|T_N[f,a]\|_{B^s_{p,r}}\leq C\|f\|_{B^s_{p,\infty}},\\ \text{$1/2<p\leq \infty$, $0<r\leq\infty$, and $1/p-1<s<\min\{1/p,1\}$},
   \end{multline}
   implies   Schauder basis properties for suitable enumerations of the Haar system. For the sake of completeness we give a sketch of this inequality which relies on the arguments in the previous section.

\begin{proof}[Proof of \eqref{T_N}]
%We proof the uniform boundedness of $T_N[\cdot,a]$ between different Besov spaces. 
%The arguments are actually quite similar to the proof of the boundedness of the $\bbE_N$ in Proposition \ref{qBanach} and Theorem \ref{expthm}. 
We may assume $\|a\|_\infty=1$. The 
%crucial
 modification of the proof of Proposition \ref{ENPN} is the fact that, due to the cancellation properties of the Haar functions participating in 
\eqref{TNdef} we can work directly with $\|T_N[f,a]\|_{B^s_{p,r}}$ (instead of $\|\bbE_N f - P_N f\|_{B^s_{p,r}}$. 
%Clearly, the described method allows for pulling out the $a_{\mu}$ on the expense of $\|a\|_{\infty}$.
 \\
\indent {\em Case 1.1.} Suppose $j+\ell, j > N$.  The estimates in \eqref{p>1}, \eqref{20_1} apply almost literally to  
$T_N[f,a]$ and yield estimates which are uniform for $\|a\|_\infty=1$.
%producing the additional factor $\|a\|_{\infty}$ on the right-hand side. 
Note, that we did not yet need any cancellation of the Haar functions. \\
\indent {\em Case 1.2.} Suppose $j+\ell > N \geq j$. We do not have to use \eqref{20_2} and work directly with 
$\|2^{js}T_N[f_{j+\ell},a]\|_p$. An analogous identity to \eqref{mu-conv} holds true with $\bbE_N(f_{j+\ell})-f_{j+\ell}$
replaced by $T_N[f_{j+\ell},a]$ due to the cancellation of the Haar functions $h_{N,\mu}$. In what follows we only have to care for a counterpart of $F_0$ since $F_1$ and $F_2$ 
do not show up. We end up with a counterpart of \eqref{F0} for $\|2^{js}T_N[f_{j+\ell},a]\|_p$.
% with an additional $\|a\|_{\infty}$ on the right-hand side. 
\\
\indent {\em Case 2.1.} Suppose $N\geq j+\ell, j$. Again, due to the cancellation of the Haar function, we obtain a version of \eqref{mu-conv} as in 
Case 1.2. The mean value theorem applied to the first factor 
in the integral gives the factor $2^{2j-2N}$, whereas the cancellation of $h_{N,\mu}$ gives $|T_N(\psi_{j+\ell,\nu})(x)| \lesssim 2^{j+\ell-N}$. We continue as in the proof of 
Proposition \ref{ENPN}. \\
\indent {\em Case 2.2.} The remaining case $j+\ell \leq N < j$ goes analogously to Case  B2.2. in the proof of Proposition \ref{ENPN}. Note, that also here the splitting in 
\eqref{final} and the subsequent consideration for the second summand on the right-hand side is not necessary. 
\end{proof}

{\bf Acknowledgment.}  The authors worked on this project while participating in the 2016 
summer program in Constructive Approximation and Harmonic Analysis at the Centre de Recerca Matem\`atica. 
They  would like to thank the organizers of the  program
for providing a pleasant and fruitful research atmosphere. The authors would also like to thank an 
anonymous referee for a careful proofreading and several valuable comments how to improve the presentation of the proofs. 
% \red{T.U. gratefully acknowledges 
%support by the German Research Foundation (DFG) and the Emmy-Noether
%programme, Ul-403/1-1.}

\end{document}